\documentclass[10pt,a4paper,draft]{article}

\usepackage{amsmath}
\usepackage{amsthm}
\usepackage{amssymb}
\usepackage{amscd}

\theoremstyle{definition}
\newtheorem{theorem}{Theorem}[section]
\newtheorem{prop}[theorem]{Proposition}
\newtheorem{lemma}[theorem]{Lemma}

\newenvironment{demo}[1]{%
  \trivlist
  \item[\hskip\labelsep
        {\bf #1.}]
}{%
\hfill\qedsymbol
  \endtrivlist
}

\newenvironment{roster}[1]{%
\begin{list}{}{%
\setlength{\topsep}{0pt}
\setlength{\itemsep}{0pt}
\setlength{\parsep}{0pt}
\settowidth{\labelwidth}{#1}
\setlength{\leftmargin}{\labelwidth}
\addtolength{\leftmargin}{\parindent}
}}{\end{list}}

\newcommand\Pf{\operatorname{Pf}}
\newcommand\Comp{\mathbb{C}}

\newcommand{\vectx}{\boldsymbol{x}}
\newcommand{\vecty}{\boldsymbol{y}}
\newcommand{\vectz}{\boldsymbol{z}}
\newcommand{\vectw}{\boldsymbol{w}}
\newcommand{\vecta}{\boldsymbol{a}}
\newcommand{\vectb}{\boldsymbol{b}}
\newcommand{\vectc}{\boldsymbol{c}}
\newcommand{\vectd}{\boldsymbol{d}}

\newcommand{\ep}{\varepsilon}

\title{
An elliptic generalization of Schur's Pfaffian identity
}

\author{
Soichi OKADA\thanks{
Graduate School of Mathematics, Nagoya University,
e-mail: okada@math.nagoya-u.ac.jp
}
}

\date{}

\begin{document}

\maketitle

\begin{abstract}
We present a Pfaffian identity involving elliptic functions,
whose rational limit gives a generalization of Schur's Pfaffian identity
 for $\Pf ( (x_j - x_i)/(x_j + x_i) )$.
This identity is regarded as a Pfaffian counterpart of
 Frobenius identity, which is an elliptic generalization of
 Cauchy's determinant identity for $\det ( 1/ (x_i + x_j) )$.
\end{abstract}

\section{
Introduction
}

Let $[x]$ denote a nonzero holomorphic function on the complex plane
 $\Comp$ in the variable $x$ satisfying the following two conditions:
\begin{roster}{(ii)}
\item[(i)]
$[x]$ is an odd function, i.e., $[-x] = - [x]$.
\item[(ii)]
$[x]$ satisfies the Riemann relation:
$$
[x+y][x-y][u+v][u-v] - [x+u][x-u][y+v][y-v] + [x+v][x-v][y+u][y-u] = 0.
$$
\end{roster}
It is known that such a function $[x]$ is obtained from one of the following
 functions by the transformation $[x] \to e^{ax^2+b}[cx]$
(see \cite[Chap.~XX, Misc. Ex. 38]{WW}):
\begin{roster}{(c)}
\item[(a)]
(elliptic case)
$[x] = \sigma(x)$ (the Weierstrass sigma function).
\item[(b)]
(trigonometric case)
$[x] = e^x - e^{-x}$.
\item[(c)]
(rational case)
$[x] = x$.
\end{roster}

G.~Frobenius \cite{F} gave the following determinant identity:
\begin{equation}
\det \left( \frac{[z + x_i + y_j]}{[z][x_i + y_j]} \right)_{1 \le i, j \le n}
 =
 \frac{\prod_{1 \le i < j \le n} [x_j - x_i][y_j - y_i] }
      {\prod_{1 \le i, j \le n} [x_i + y_j] }
 \cdot
 \frac{[z + \sum_{i=1}^n x_i + \sum_{j=1}^n y_j]}
      {[z]}.
\label{frobenius}
\end{equation}
If we take the limit $z \to \infty$ in the rational case of
 this identity (\ref{frobenius}), 
we obtain the Cauchy's determinant identity \cite{C}:
\begin{equation}
\det \left( \frac{1}{x_i + y_j} \right)_{1 \le i, j \le n}
 =
\frac{\prod_{1 \le i < j \le n} (x_j - x_i) (y_j - y_i) }
     {\prod_{1 \le i, j \le n} (x_i + y_j) }.
\label{cauchy}
\end{equation}
Hence the identity (\ref{frobenius}) can be regarded as
 an elliptic generalization of (\ref{cauchy}).

Among identities for Pfaffians, the Schur's Pfaffian identity
 \cite[p.~225]{S}
\begin{equation}
\Pf \left( \frac{x_i - x_i}{x_i + x_i} \right)_{1 \le i, j \le 2n}
 =
\prod_{1 \le i < j \le 2n} \frac{x_i - x_i}{x_j + x_i}
\label{schur}
\end{equation}
plays the same role in the theory of Schur's $Q$-functions as
 the Cauchy's determinant identity in the theory of Schur functions.
Our main result is the following elliptic generalization of
 the Schur's Pfaffian identity.

\begin{theorem} \label{thm:main}
For complex variables $x_1, \cdots, x_{2n}$, $z$ and $w$, we have
\begin{multline}
\Pf \left(
 \frac{ [x_j - x_i] }{ [x_j + x_i] }
 \cdot \frac{ [z + x_i + x_j] }{ [z] }
 \cdot \frac{ [w + x_i + x_j] }{ [w] }
\right)_{1 \le i, j \le 2n}
\\
 =
 \prod_{1 \le i < j \le 2n} \frac{ [x_j - x_i] }{ [x_j + x_i] }
 \cdot \frac{ [z + \sum_{i=1}^{2n} x_i] }{ [z] }
 \cdot \frac{ [w + \sum_{i=1}^{2n} x_i] }{ [w] }.
\label{main}
\end{multline}
\end{theorem}

If we take the limit $z \to \infty$, $w \to \infty$ in the rational case
 of (\ref{main}),
 then we recover the Schur's Pfaffian identity (\ref{schur}).

Recently, the Frobenius' identity (\ref{frobenius}) plays a key role
 in proving transformations of elliptic hypergeometric series.
See \cite{KN}, \cite{R}.
We expect that the identity (\ref{main}) provides another tool in the theory
 of special functions.

This paper is organized as follows.
In Section~2, we give a proof of Theorem~\ref{thm:main}.
In Section~3, we present a generalization of the rational case of
 (\ref{frobenius}) and (\ref{main}), in terms of Schur functions
 corresponding to the staircase partitions.
And we give another proof of the trigonometric case of 
 (\ref{frobenius}) and (\ref{main}) in Section~4.

\section{
Proof of Theorem~\ref{thm:main}
}

In this section, we give a proof of Theorem~\ref{thm:main}.
The key is to to use the following Pfaffian version of the Desnanot--Jacobi
 formula.
See \cite[(1.1)]{Kn} and \cite[Theorem~2.6]{IW} for a proof and
 related formulae.

\begin{lemma} \label{pf-dodgson}
Let $A$ be the skew-symmetric matrix and denote by $A^{i_1, \cdots, i_k}$
 the skew-symmetric matrix obtained by removing $i_1$th, $\cdots$, $i_k$th
 rows and columns.
Then we have
$$
\Pf A^{12} \cdot \Pf A^{34} - \Pf A^{13} \cdot \Pf A^{24}
 + \Pf A^{14} \cdot \Pf A^{23}
 = \Pf A \cdot \Pf A^{1234}.
$$
\end{lemma}

\begin{demo}{Proof of Theorem~\ref{thm:main}}
We proceed by induction on $n$.
By applying Lemma~\ref{pf-dodgson} to the skew-symmetric matrix
 on the left hand side of (\ref{main}),
and using the induction hypothesis, we see that it is enough to show
\begin{align}
&
[b-a][d-c][c+a][d+a][c+b][d+b]
\notag
\\
&\quad\times
[z+a+b+s][z+c+d+s][w+a+b+s][w+c+d+s]
\notag
\\
&
-
[c-a][d-b][b+a][d+a][c+b][d+c]
\notag
\\
&\quad\times
[z+a+c+s][z+b+d+s][w+a+c+s][w+b+d+s]
\notag
\\
&
+
[d-a][c-b][b+a][c+a][d+b][d+c]
\notag
\\
&\quad\times
[z+a+d+s][z+b+c+s][w+a+d+s][w+b+c+s]
\notag
\\
&
=
[b-a][c-a][d-a][c-b][d-b][d-c]
\notag
\\
&\quad\times
[z+a+b+c+d+s][z+s][w+a+b+c+d+s][w+s],
\label{key}
\end{align}
where $a = x_1$, $b = x_2$, $c = x_3$, $d = x_4$ and $s = \sum_{j=5}^{2n} x_j$.
Note that the case of $s=0$ is exactly the identity (\ref{main}) with $n=2$.

By replacing $z+s$ and $w+s$ by $z$ and $w$ respectively,
 we may assume $s=0$.
By applying the Riemann relation with $(x,y,u,v) = (z/2+a,z/2+b,z/2+c,z/2+d)$,
 we have
$$
[b-a][z+a+b][d-c][z+c+d]
= [c-a][z+a+c][d-b][z+b+d]
 - [d-a][z+a+d][c-b][z+b+c].
$$
Hence the left hand side of (\ref{key}) is equal to
\begin{align}
&
[c-a][z+a+c][d-b][z+b+d][d+a][c+b]
\notag
\\
&\quad
\times 
\left( [c+a][d+b][w+a+b][w+c+d] - [b+a][d+c][w+a+c][w+b+d] \right)
\notag
\\
&
-
[d-a][z+a+d][c-b][z+b+c][c+a][d+b]
\notag
\\
&\quad
\times
\left( [d+a][c+b][w+a+b][w+c+d] - [b+a][d+c][w+a+d][w+b+c] \right).
\label{step1}
\end{align}
Again, by using the Riemann relation with
$$
x = \frac{w+2a+b+c}{2},
\quad
y = \frac{w+b-c}{2},
\quad
u = \frac{w+b+c+2d}{2},
\quad
v = \frac{w-b+c}{2},
$$
we have
\begin{multline*}
[c+a][d+b][w+a+b][w+c+d] - [b+a][d+c][w+a+c][w+b+d]
\\
=
[d-a][c-b][w][w+a+b+c+d].
\end{multline*}
By exchanging $c$ and $d$, we have
\begin{multline*}
[d+a][c+b][w+a+b][w+c+d] - [b+a][d+c][w+a+d][w+b+c]
\\
=
[c-a][d-b][w][w+a+b+c+d].
\end{multline*}
Hence we see that (\ref{step1}) is equal to
\begin{multline}
[c-a][d-a][c-b][d-b][w][w+a+b+c+d]
\\
\times
\left(
[d+a][c+b][z+a+c][z+b+d]
-
[c+a][d+b][z+a+d][z+b+c]
\right).
\label{step2}
\end{multline}
Finally we use the Riemann relation with
$$
x = \frac{z+2a+c+d}{2},
\quad
y = \frac{z+c-d}{2},
\quad
u = \frac{z+2b+c+d}{2},
\quad
v = \frac{z-c+d}{2}.
$$
Then we have
\begin{multline*}
[d+a][c+b][z+a+c][z+b+d]
-
[c+a][d+b][z+a+d][z+b+c]
\\
=
[b-a][d-c][z][z+a+b+c+d],
\end{multline*}
and we see that (\ref{step2}) is equal to the right hand side of (\ref{key}).
This completes the proof of Theorem~\ref{thm:main}.
\end{demo}

\section{
Rational case
}

In this section, we generalize the rational case of
 (\ref{frobenius}) and (\ref{main}) to the identities
 involving Schur functions.

If we put $[x] = x$, then the identities (\ref{frobenius}) and
 (\ref{main}) reduce to
\begin{align}
&
\det \left(
 \frac{1}{x_i+y_j} \cdot \frac{z+x_i+y_j}{z}
\right)_{1 \le i, j \le n}
\notag
\\
&\quad=
\frac{ \prod_{1 \le i < j \le n} (x_j - x_i) (y_j - y_i) }
     { \prod_{1 \le i, j \le n} (x_i + y_j) }
\cdot
\frac{ z + \sum_{i=1}^n x_i + \sum_{j=1}^n y_j }{z},
\label{rational1}
\\
&
\Pf \left(
 \frac{ x_j - x_i }{ x_j + x_i }
 \cdot \frac{ z + x_i + x_j }{ z }
 \cdot \frac{ w + x_i + x_j }{ w }
\right)_{1 \le i, j \le 2n}
\notag
\\
&\quad=
 \prod_{1 \le i < j \le 2n} \frac{ x_j - x_i }{ x_j + x_i }
 \cdot \frac{ z + \sum_{i=1}^{2n} x_i }{ z }
 \cdot \frac{ w + \sum_{i=1}^{2n} x_i }{ w },
\label{rational2}
\end{align}
respectively.
These identities can be generalized as follows:

\begin{theorem} \label{thm:rational}
Let $s_\lambda$ be the Schur function corresponding to a partition
 $\lambda$.
For a nonnegative integer $r$, let $\delta(r) =(r,r-1,\cdots, 2,1)$ denote
 the staircase partition.
\begin{roster}{(1)}
\item[(1)]
Let $k \le m$ be integers.
For vectors of variables $\vectx = (x_1, \cdots, x_n)$, $\vecty
 = (y_1, \cdots, y_n)$ and $\vectz = (z_1, \cdots, z_m)$, we have
\begin{equation}
\det \left(
 \frac{ 1 }{ x_i + y_j }
 \cdot
 \frac{ s_{\delta(k)}(x_i, y_j, \vectz) }
      { s_{\delta(k)}(\vectz) }
\right)_{1 \le i, j \le n}
 =
\frac{ \prod_{1 \le i < j \le n} (x_j - x_i) (y_j - y_i) }
     { \prod_{1 \le i, j \le n} (x_i + y_j) }
\cdot
\frac{ s_{\delta(k)}(\vectx, \vecty, \vectz) }{ s_{\delta(k)}(\vectz) }.
\label{rational3}
\end{equation}
\item[(2)]
Let $k \le m$ and $l \le m'$ be integers.
For vectors of variables $\vectx = (x_1, \cdots, x_{2n})$,
 $\vectz = (z_1, \cdots, z_m)$ and $\vectw = (w_1, \cdots, w_{m'})$,
 we have
\begin{multline}
\Pf \left(
 \frac{ x_j - x_i }{ x_j + x_i }
 \cdot
 \frac{ s_{\delta(k)}(x_i,x_j,\vectz) }{ s_{\delta(k)}(\vectz) }
 \cdot
 \frac{ s_{\delta(l)}(x_i,x_j,\vectw) }{ s_{\delta(l)}(\vectw) }
\right)_{1 \le i, j \le 2n}
\\
 =
 \prod_{1 \le i < j \le 2n} \frac{ x_j - x_i }{ x_j + x_i }
 \cdot
 \frac{ s_{\delta(k)}(\vectx,\vectz) }{ s_{\delta(k)}(\vectz) }
 \cdot
 \frac{ s_{\delta(l)}(\vectx,\vectw) }{ s_{\delta(l)}(\vectw) }.
\label{rational4}
\end{multline}
\end{roster}
\end{theorem}

The identity (\ref{rational3}) reduces to (\ref{cauchy}) and (\ref{rational1})
 in the case of $k=0$ and $k=1$ respectively,
 while (\ref{rational4}) gives (\ref{schur}) and (\ref{rational2})
 in the case of $k=l=0$ and $k=l=1$ respectively.

We derive Theorem~\ref{thm:rational} from much more general identities
 involving generalized Vandermonde determinants, which were conjectured
 by the author and proven in \cite{IOTZ}.
Let $\vectx = (x_1, \cdots, x_n)$ and $\vecta = (a_1, \cdots, a_n)$
 be two vectors of variables of length $n$.
For nonnegative integers $p$ and $q$ with $p+q = n$,
 we define a generalized Vandermonde matrix $V^{p,q}(\vectx ; \vecta)$
 to be the $n \times n$ matrix with $i$th row
$$
(1, x_i, \cdots, x_i^{p-1}, a_i, a_i x_i, \cdots, a_i x_i^{q-1}).
$$
Then we have

\begin{theorem} \label{thm:general}
 (Ishikawa--Okada--Tagawa--Zeng \cite{IOTZ})
\begin{roster}{(d)}
\item[(a)]
Let $n$ be a positive integer and let $p$ and $q$ be nonnegative integers.
For six vectors of variables
\begin{gather*}
\vectx = (x_1, \cdots, x_n),\ 
\vecty = (y_1, \cdots, y_n),\ 
\vecta = (a_1, \cdots, a_n),\ 
\vectb = (b_1, \cdots, b_n),
\\
\vectz = (z_1, \cdots, z_{p+q}),\ 
\vectc = (c_1, \cdots, c_{p+q}),
\end{gather*}
we have
\begin{multline}
\det \left(
 \frac{ \det V^{p+1,q+1}(x_i,y_j,\vectz ; a_i,b_j,\vectc) }
      {y_j - x_i}
 \right)_{1 \le i, j \le n}
\\
=
\frac{ (-1)^{n(n-1)/2} }{ \prod_{i,j=1}^n (y_j - x_i) }
 \det V^{p,q}(\vectz ; \vectc)^{n-1}
 \det V^{n+p,n+q}(\vectx,\vecty,\vectz ; \vecta,\vectb,\vectc).
\label{general1}
\end{multline}
\item[(b)]
Let $n$ be a positive integer and
 let $p$, $q$, $r$, $s$ be nonnegative integers.
For seven vectors of variables
\begin{gather*}
\vectx = (x_1, \cdots, x_{2n}),\ 
\vecta = (a_1, \cdots, a_{2n}),\ 
\vectb = (b_1, \cdots, b_{2n}),
\\
\vectz = (z_1, \cdots, z_{p+q}),\ 
\vectc = (c_1, \cdots, c_{p+q}),
\\
\vectw = (w_1, \cdots, w_{r+s}),\ 
\vectd = (d_1, \cdots, d_{r+s}),
\end{gather*}
we have
\begin{multline}
\Pf \left(
 \frac{ \det V^{p+1,q+1}(x_i,x_j,\vectz ; a_i,a_j,\vectc)
        \det V^{r+1,s+1}(x_i,x_j,\vectw ; b_i,b_j,\vectd) }
      { x_j - x_i }
\right)_{1 \le i, j \le 2n}
\\
=
\frac{1}{\prod_{1 \le i < j \le 2n}(x_j - x_i)}
 \det V^{p,q}(\vectz ; \vectc)^{n-1}
 \det V^{r,s}(\vectw ; \vectd)^{n-1}
\\
\times
 \det V^{n+p,n+q}(\vectx,\vectz ; \vecta,\vectc)
 \det V^{n+r,n+s}(\vectx,\vectw ; \vectb,\vectd).
\label{general2}
\end{multline}
\end{roster}
\end{theorem}

\begin{demo}{Proof of Theorem~\ref{thm:rational}}
In (\ref{general1}) (resp. (\ref{general2})), we specialize
\begin{equation}
\begin{gathered}
x_i \to x_i^2,
\quad
y_i \to y_i^2,
\quad
z_i \to z_i^2,
\quad
a_i \to x_i,
\quad
b_i \to y_i,
\quad
c_i \to z_i,
\\
(\text{resp.}
\quad
x_i \to x_i^2,
\quad
z_i \to z_i^2,
\quad
w_i \to w_i^2,
\quad
a_i \to x_i,
\quad
b_i \to x_i,
\quad
c_i \to z_i,
\quad
d_i \to w_i).
\label{specialization}
\end{gathered}
\end{equation}
It follows from the bi-determinant definition of Schur functions that
$$
\det V^{p,q}(x_1^2, \cdots, x_{p+q}^2 ; x_1, \cdots, x_{p+q})
 =
\ep_{p,q} s_{\delta_{p,q}}(x_1, \cdots, x_{p+q})
\prod_{1 \le i < j \le p+q} (x_j - x_i),
$$
where the partition $\delta_{p,q}$ and the signature $\ep_{p,q}$ are
 defined by
$$
\delta_{p,q}
 = \begin{cases}
 \delta(p-q-1) &\text{if $p > q$,} \\
 \delta(q-p) &\text{if $p \le q$,}
\end{cases}
\quad
\ep_{p,q}
 = \begin{cases}
 (-1)^{q(2p-q-1)/2} &\text{if $p > q$,} \\
 (-1)^{p(p-1)/2} &\text{if $p \le q$.}
\end{cases}
$$
Hence, under the specialization (\ref{specialization}), the identities
 (\ref{general1}) and (\ref{general2}) in Theorem~\ref{thm:general} give
\begin{align*}
&
\ep_{p+1,q+1}^n
\det \left(
 \frac{ s_{\delta_{p+1,q+1}}(x_i, y_j, \vectz) }
      { x_i + y_j }
 \right)_{1 \le i, j \le n}
\\
&\quad=
 \frac{ (-1)^{n(n-1)/2} \ep_{p,q}^{n-1} \ep_{n+p,n+q} }
      { \prod_{i,j=1}^n (x_i + y_j) }
 \cdot s_{\delta_{p,q}}(\vectz)^{n-1}
 \cdot s_{\delta_{n+p,n+q}}(\vectx,\vecty,\vectz),
\\
&
\ep_{p+1,q+1}^n \ep_{r+1,s+1}^n
\Pf \left(
 \frac{x_j - x_i}{x_j + x_i}
 \cdot s_{\delta_{p,q}}(x_i, x_j, \vectz)
 \cdot s_{\delta_{r,s}}(x_i, x_j, \vectw)
 \right)_{1 \le i, j \le 2n}
\\
&\quad=
\ep_{p,q}^{n-1} \ep_{r,s}^{n-1} \ep_{p+n,q+n} \ep_{r+n,s+n}
\prod_{1 \le i < j \le 2n} \frac{x_j - x_i}{x_j + x_i}
\\
&\quad\quad\times
s_{\delta_{p,q}}(\vectz)^{n-1}
\cdot s_{\delta_{r,s}}(\vectw)^{n-1}
\cdot s_{\delta_{n+p,n+q}}(\vectx,\vectz)
\cdot s_{\delta_{n+r,n+s}}(\vectx,\vectw).
\end{align*}
Here a direct computation shows
$$
(-1)^{n(n-1)/2} \ep_{p,q}^{n-1} \ep_{n+p,n+q}
 = \ep_{p+1,q+1}^n.
$$

Now we can complete the proof by choosing $(p,q)$ and $(p,q,r,s)$ as follows.
In the determinant case, given integers $m$ and $k$, we define $p$ and $q$ by
$$
(p,q) = \begin{cases}
 \left( \dfrac{m-k}{2}, \dfrac{m+k}{2} \right)
 &\text{if $m$ and $k$ have the same parity,} \\
 \left( \dfrac{m+k+1}{2}, \dfrac{m-k-1}{2} \right)
 &\text{otherwise.}
\end{cases}
$$
In the Pfaffian case, for integers $m$, $m'$, $k$ and $l$, we put
\begin{gather*}
(p,q) = \begin{cases}
 \left( \dfrac{m-k}{2}, \dfrac{m+k}{2} \right)
 &\text{if $m$ and $k$ have the same parity,} \\
 \left( \dfrac{m+k+1}{2}, \dfrac{m-k-1}{2} \right)
 &\text{otherwise,}
\end{cases}
\\
(r,s) = \begin{cases}
 \left( \dfrac{m'-l}{2}, \dfrac{m'+l}{2} \right)
 &\text{if $m'$ and $l$ have the same parity,} \\
 \left( \dfrac{m'+l+1}{2}, \dfrac{m'-l-1}{2} \right)
 &\text{otherwise.}
\end{cases}
\end{gather*}
\end{demo}

\section{
Trigonometric case
}

In this section, we give another proof of the trigonometric case
 of (\ref{frobenius}) and (\ref{main}) by using determinant
 and Pfaffian identities similar to (\ref{general1}) and (\ref{general2}).

By replacing $[x]$ by $x^{1/2}-x^{-1/2}$, we can see that the trigonometric
 cases of (\ref{frobenius}) and (\ref{main}) are equivalent to
\begin{align}
&
\det \left( \frac{1 - z x_i y_j}{1 - x_i y_j} \right)_{1 \le i, j \le n}
\notag
\\
&\quad=
(1 - z)^{n-1} \left( 1 - z \prod_{i=1}^n x_i \prod_{j=1}^n y_j \right)
\cdot
\frac{ \prod_{1 \le i < j \le n} (x_j - x_i)(y_j - y_i) }
     { \prod_{1 \le i, j \le n} (1 - x_i y_j) },
\label{trigonometric1}
\\
&
\Pf \left(
 \frac{(x_j - x_i)(1 - z x_i x_j)(1 - w x_i x_j)}{1 - x_i x_j}
\right)_{1 \le i, j \le 2n}
\notag
\\
&\quad=
(1-z)^{n-1} (1-w)^{n-1}
\left( 1 - z \prod_{i=1}^{2n} x_i \right)
\left( 1 - w \prod_{i=1}^{2n} x_i \right)
\cdot
\prod_{1 \le i < j \le 2n} \frac{x_i - x_i}{1 - x_i x_j},
\label{trigonometric2}
\end{align}
respectively.
Here we derive these identities from the following Proposition.

\begin{prop} (Okada \cite{O})
\begin{roster}{(2)}
\item[(1)]
For vectors of variables $\vectx$, $\vecty$, $\vecta$ and $\vectb$
 of length $n$, we have
\begin{align}
&
\det \left( \frac{ 1 - a_i b_j }{ 1 - x_i y_j } \right)_{1 \le i, j \le n}
\notag
\\
&\quad
 = \frac{ (-1)^{n(n+1)/2} }
        { \prod_{1 \le i, j \le n} (1 - x_i y_j) }
\det
\begin{pmatrix}
a_1 x_1^{n-1} & a_1 x_1^{n-2} & \cdots & a_1
 & x_1^{n-1} & x_1^{n-2} & \cdots & 1 \\
\vdots & \vdots & & \vdots
 & \vdots & \vdots & & \vdots \\
a_n x_n^{n-1} & a_n x_n^{n-2} & \cdots & a_n
 & x_n^{n-1} & x_n^{n-2} & \cdots & 1 \\
1 & y_1 & \cdots & y_1^{n-1}
 & b_1 & b_1 y_1 & \cdots & b_1 y_1^{n-1} \\
\vdots & \vdots & & \vdots
 & \vdots & \vdots & & \vdots \\
1 & y_n & \cdots & y_n^{n-1}
 & b_n & b_n y_n & \cdots & b_n y_n^{n-1}
\end{pmatrix}.
\label{general3}
\end{align}
\item[(2)]
For vectors of variables $\vectx$, $\vecta$ and $\vectb$
 of length $2n$, we have
\begin{multline}
\Pf \left(
 \frac{\det W^2(x_i, x_j ; a_i, a_j) \det W^2(x_i, x_j ; b_i, b_j )}
      {(x_j - x_i)(1 - x_i x_j)}
\right)_{1 \le i, j \le 2n}
\\
 =
\frac{1}{\prod_{1 \le i < j \le 2n} (x_j - x_i)(1 - x_i x_j)}
 \det W^{2n}(\vectx;\vecta) \det W^{2n}(\vectx;\vectb),
\label{general4}
\end{multline}
where $W^n(\vectx ; \vecta)$ is the $n \times n$ matrix with $i$th row
$$
(1+a_i x_i^{n-1}, x_i + a_i x_i^{n-2}, \cdots, x_i^{n-1} + a_i).
$$
\end{roster}
\end{prop}

\begin{demo}{Proof}
The first identity (\ref{general1}) is obtained by replacing
 $a_i$ and $x_i$ by $1/a_i$ and $1/x_i$ respectively
 in the identity (\ref{general1}) with $p=q=0$.
 (This case was first given in \cite[Theorem~4.2]{O}.)
The second identity (\ref{general2}) is given in \cite[Theorem~4.4]{O}.
\end{demo}

Now the trigonometric case (\ref{trigonometric1})
 (resp. (\ref{trigonometric2}))
 easily follows from (\ref{general3}) (resp. (\ref{general2}))
 by substituting $a_i = z^{1/2} x_i$ and $b_i = z^{1/2} y_i$
 (resp. $a_i = z x_i$ and $b_i = w x_i$),
 and by using the determinant evaluation in Lemma~\ref{lem:det1}
 (resp. \ref{lem:det2}).

\begin{lemma} \label{lem:det1}
\begin{multline}
\det
\begin{pmatrix}
t x_1^n & t x_1^{n-1} & \cdots & t x_1
 & x_1^{n-1} & x_1^{n-2} & \cdots & 1 \\
\vdots & \vdots & & \vdots
 & \vdots & \vdots & & \vdots \\
t x_n^n & t x_n^{n-1} & \cdots & t x_n
 & x_n^{n-1} & x_n^{n-2} & \cdots & 1 \\
1 & y_1 & \cdots & y_1^{n-1}
 & t y_1 & t y_1^2 & \cdots & t y_1^n \\
\vdots & \vdots & & \vdots
 & \vdots & \vdots & & \vdots \\
1 & y_n & \cdots & y_n^{n-1}
 & t y_n & t y_n^2 & \cdots & t y_n^n
\end{pmatrix}
\\
= (-1)^{n(n+1)/2} (1-t^2)^{n-1}
 \left( 1 - t^2 \prod_{i=1}^n x_i \prod_{i=1}^n y_i \right)
 \prod_{1 \le i < j \le n} (x_j - x_i)(y_j - y_i).
\label{det1}
\end{multline}
\end{lemma}

\begin{demo}{Proof}
Let $D(\vectx, \vecty ; t)$ denote the determinant of the left hand side.
Since $D(\vectx,\vecty;t)$ is divisible by
 $\Delta(\vectx) \Delta(\vecty) =
 \prod_{1 \le i < j \le n} (x_j - x_i) (y_j - y_i)$ and it has degree $n$
 with respect to $x_1$, we can write
\begin{equation}
D(\vectx, \vecty ; t)
 = (c_0 + c_1 x_1) \prod_{1 \le i < j \le n} (x_j - x_i)(y_j - y_i),
\label{7}
\end{equation}
where $c_0$ and $c_1$ are polynomials in $x_2, \cdots, x_n, y_1, \cdots, y_n,
 t$.

By performing elementary transformations and by using the Vandermonde
 determinant, we see that the constant term and the coefficient of $x_1^n$
 in $D(\vectx, \vecty ; t)$ are given by
\begin{gather*}
(1-t^2)^{n-1} (-1)^{(n-1)(n-2)/2}
 \prod_{2 \le i < j \le n}(x_j - x_i)
 \prod_{1 \le i < j \le n}(y_j - y_i),
\\
t (1-t^2)^{n-1} (-1)^{n(n-1)/2} \prod_{i=2}^n x_i
 \prod_{2 \le i < j \le n} (x_j - x_i)
 \prod_{1 \le i < j \le n} (y_j - y_i),
\end{gather*}
respectively.
These computation yield the coefficients $c_0$ and $c_1$:
$$
c_0 = (-1)^{n(n+1)/2} (1-t^2)^{n-1},
\quad
c_1 = - (-1)^{n(n+1)/2} (1-t^2)^{n-1} t^2 \prod_{i=2}^n x_i \prod_{i=1}^n y_i,
$$
which complete the proof of (\ref{det1}).
\end{demo}

\begin{lemma} \label{lem:det2}
Let $d_n(t)$ be the polynomial given by
$$
d_n(t) = \begin{cases}
 (1-t)^m (1+t)^m &\text{if $n = 2m$ is even,} \\
 (1-t)^m (1+t)^{m+1} &\text{if $n = 2m+1$ is odd.}
\end{cases}
$$
Then we have
\begin{align}
\det W^n(x_1, \cdots, x_n ; t, \cdots, t)
 &=
d_n(t) \prod_{1 \le i < j \le n} (x_j - x_i),
\\
\det W^n(x_1, \cdots, x_n ; tx_1, \cdots, tx_n)
 &=
d_{n-1}(t)
 \left( 1 - (-1)^n t \prod_{i=1}^n x_i \right)
 \prod_{1 \le i < j \le n} (x_j - x_i).
\label{det2}
\end{align}
\end{lemma}

\begin{demo}{Proof}
The idea of the proof is the same as Lemma~\ref{lem:det1},
 so we leave it to the reader.
\end{demo}


\begin{thebibliography}{99}

\bibitem{C}
A.~L.~Cauchy,
M\'emoire sur les fonctions altern\'ees et sur les sommes altern\'ees,
Exercices Anal. et Phys. Math. {\bf 2} (1841), 151--159.

\bibitem{F}
G.~Frobenius,
\"Uber die elliptischen Funktionen zweiter Art,
J. f\"ur die reine und ungew. Math. {\bf 93} (1882), 53--68.

\bibitem{IOTZ}
M.~Ishikawa, S.~Okada, H.~Tagawa and J.~Zeng, 
Generalizations of Cauchy's determinant and Schur's Pfaffian,
{\tt arXiv:math.CO/0411280}.

\bibitem{IW}
M.~Ishikawa and M.~Wakayama,
Applications of the minor summation formula III :
 Pl\"ucker relations, lattice paths and Pfaffians,
{\tt arXiv:math.CO/0312358}.

\bibitem{KN}
Y.~Kajiwara and M.~Noumi,
Multiple elliptic hypergeometric series : an approach from
 the Cauchy determinant,
Indag. Math. (N.S.) {\bf 14} (2003), 395--421.

\bibitem{Kn}
D.~Knuth,
Overlapping Pfaffians,
Electron. J. Combin. {\bf 3} (2) (``The Foata Festschrift'') (1996),
 151--163.

\bibitem{O}
S.~Okada,
Applications of minor-summation formulas to rectangular-shaped representations
 of classical groups,
J. Alg. {\bf  205} (1998), 337--367.

\bibitem{R}
E.~M.~Rains,
Transformations of elliptic hypergeometric integrals,
{\tt arXiv:math.QA/0309252}.

\bibitem{S}
I.~Schur,
\"Uber die Darstellung der symmetrischen und der alternirenden Gruppe
 durch gebrochene lineare Substitutuionen,
J. Reine Angew. Math. {\bf 139} (1911), 155--250.

\bibitem{WW}
E.~T.~Whittaker and G.~N.~Watson,
A Course of Modern Analysis (Fourth edition),
Cambridge Univ. Press, 1927 (reprinted in 1999).

\end{thebibliography}
\end{document}